\documentclass[11pt]{article}

\usepackage{epsfig}
\usepackage{amsfonts,amssymb,latexsym}

\newtheorem{theorem}{Theorem}[section]

\newtheorem{defi}{Definition}[section]
\newtheorem{lemma}[theorem]{Lemma}

\def\binom#1#2{{#1}\choose{#2}}

\def\slfrac#1#2{\hbox{\kern.1em %
 \raise.5ex\hbox{\the\scriptfont0 #1}\kern-.11em %
 /\kern-.15em\lower.25ex\hbox{\the\scriptfont0 #2}}}

\newcommand{\ep}{\epsilon}
\newcommand{\hsp}{\hspace*{\parindent}}

\newcommand{\beq}{\begin{eqnarray}}
\newcommand{\eeq}{\end{eqnarray}}
\newcommand{\beql}[1]{\begin{eqnarray}\label{#1}}
\newcommand{\beqs}{\begin{eqnarray*}}
\newcommand{\eeqs}{\end{eqnarray*}}
\newcommand{\eqn}[1]{(\ref{#1})}

\newcommand{\lf}{\lfloor}
\newcommand{\rf}{\rfloor}

\newcommand{\cc}{{\mathbb C}}

\newcommand{\pp}{{\mathbb P}}
\newcommand{\rr}{{\mathbb R}}
\newcommand{\zz}{{\mathbb Z}}
\newcommand{\qq}{{\mathbb Q}}

\newcommand{\TT}{{\mathbb T}}

\newcommand{\sS}{{\mathcal S}}

\newcommand{\sG}{{\mathcal G}}

\newcommand{\sO}{{\mathcal O}}

\newcommand{\sW}{{\mathcal W}}

\newcommand{\bsq}{\vrule height .9ex width .8ex depth -.1ex}

\makeatletter
\def\@sect#1#2#3#4#5#6[#7]#8{\ifnum #2>\c@secnumdepth
     \def\@svsec{}\else
     \refstepcounter{#1}\edef\@svsec{\csname the#1\endcsname.\hskip .75em }\fi
     \@tempskipa #5\relax
      \ifdim \@tempskipa>\z@
        \begingroup #6\relax
          \@hangfrom{\hskip #3\relax\@svsec}{\interlinepenalty \@M #8\par}%
        \endgroup
       \csname #1mark\endcsname{#7}\addcontentsline
         {toc}{#1}{\ifnum #2>\c@secnumdepth \else
                      \protect\numberline{\csname the#1\endcsname}\fi
                    #7}\else
        \def\@svsechd{#6\hskip #3\@svsec #8\csname #1mark\endcsname
                      {#7}\addcontentsline
                           {toc}{#1}{\ifnum #2>\c@secnumdepth \else
                             \protect\numberline{\csname the#1\endcsname}\fi
                       #7}}\fi
     \@xsect{#5}}
\makeatother

\catcode`\@=11
\renewcommand{\section}{
        \setcounter{equation}{0}
        \@startsection {section}{1}{\z@}{-3.5ex plus -1ex minus
        -.2ex}{2.3ex plus .2ex}{\large\bf}
        }
\catcode`@=12

\begin{document}
\begin{center}
{\Large {\bf On the Normality of Arithmetical Constants }}\\

\vspace{1.5\baselineskip}
{\em Jeffrey C. Lagarias} \\
\vspace*{1.5\baselineskip}
(February 16, 2001 version) 
\vspace*{1.5\baselineskip}\end{center}
\noindent{\bf ABSTRACT.}
Bailey and Crandall~\cite{BC00} recently formulated ``Hypothesis $A$'',
a general principle to explain the (conjectured) normality of 
the binary expansion of 
constants like $\pi$ and other related numbers, or more generally
the base $b$ expansion of such constants for an integer $b \geq 2$.
This paper 
shows that a basic mechanism underlying their principle,
which is a relation between single orbits of two discrete 
dynamical systems,
holds for a very general class of representations of numbers.
This general class includes numbers for which 
the conclusion of  ``Hypothesis $A$'' is not true.
The paper also
relates the particular class of arithmetical constants treated by
Bailey and Crandall
to special values of $G$-functions, and points out
an analogy of ``Hypothesis $A$'' with Furstenberg's conjecture on invariant
measures.

\vspace*{1.5\baselineskip}


\section{Introduction}
\hsp
Much is known about the irrationality and
transcendence of classical arithmetical constants
such as $\pi, e,$ and $ \zeta(n)$ for $n \geq 2.$
There are general methods 
which in many cases establish
irrationality or transcendence of such numbers.
In contrast, almost
nothing is known concerning the question of whether
 arithmetical
constants are normal numbers to a fixed base, say $b = 2$.
It is unknown whether any algebraic
number is normal to any integer base $b \geq 2.$ 
Even very
weak assertions in the direction of normality are unresolved.
For example, it is not known whether $\sqrt{2}$ has arbitrarily
long blocks of zeros appearing in its binary expansion,
i.e. whether $\liminf_{n \to \infty}\{ \{2^n \sqrt{2}\}\} = 0$.

Recently
Bailey and Crandall~\cite{BC00} 
formulated  ``Hypothesis $A$'',
which provides a hypothetical general principle to explain the (conjectured)
normality to base $2$ of certain arithmetical
constants such as $\pi$ and $\log{2}.$ 

\noindent\paragraph{ Hypothesis A.}
{\em 
 Given a positive integer $b \geq 2$ and a rational
function $R(x) = \frac {p(x)}{q(x)} \in \qq(x),$
such that 
$ \deg (p(x)) < \deg (q(x))$, and  $q(x)$ has no nonnegative
integer roots, 
define $\theta = \sum_{n=0}^\infty \frac{p(n)}{q(n)} b^{-n}.$
If $y_0 = 0$ and
\beql{101a} 
y_{n + 1} = b y_{n} + \frac{p(n)}{q(n)} ~~(\bmod~ 1),
\eeq
then the sequence $\{ y_n :~n \geq 1\}$ either has finitely many
limit points or is uniformly distributed $(\bmod~ 1).$
} \\

This hypothesis
concerns the behavior of a particular orbit of the
discrete dynamical system  \eqn{101a}.
Assuming  ``Hypothesis  A,'' Bailey and Crandall deduced 
that the number $\theta$ either is rational or else is a normal
number to base $b$; these 
correspond to the
two possible behaviors of the sequence  $\{ y_n :~n \geq 1\}$
allowed by ``Hypothesis A'',  see Theorem~\ref{th41} below.
Proving ``Hypothesis A''
appears intractable, 
but it seems useful in collecting a number of 
disparate phenomona together under a single principle.
A formulation in terms of dynamical systems is natural,
because the property of normality is itself expressable
in terms of dynamics of an orbit of
another dynamical system, the $b$-transformation, see \S2.
The basic mechanism rendering ``Hypothesis A'' useful is
a relation between particular orbits of these 
two different dynamical systems.

This paper provides some complements to the results 
of Bailey and Crandall. It 
shows that the relation between particular orbits of two
discrete dynamical systems underlying ``Hypothesis A''
is valid very generally,
in that it applies to   
expansions of real numbers of the form,
$\theta = \sum_{n=1}^\infty \epsilon_n b^{-n},$
with $\epsilon_n$ arbitrary real numbers with
$\epsilon_n \to 0$ as $n \to \infty,$ 
see Theorem~\ref{th31a}.  
Every real number has such an expansion, and ``Hypothesis A''
is not true in such generality. Thus in order to be valid
``Hypothesis A'' must be restricted to  
apply only to expansions of some special form. 
Bailey and Crandall do this, formulating  ``Hypothesis A'' 
only for a countable
class of arithmetical constants which in the sequel we call
{\em BBP-numbers.}
It does not seem clear what should be the ``optimal''
class of arithmetical constants
for which ``Hypothesis A'' might
be valid. The remainder of the paper discusses
various mathematical topics relevant to this issue.
We relate $BBP$ numbers to the
theory of $G$-functions and characterize
the subclass of $BBP$-numbers which are ``special values''
of $G$-functions. We also compare ``Hypothesis A'' 
to  a conjecture of Furstenberg
in ergodic theory, and this suggests some further
questions to pursue. 

We now summarize the contents of the paper in more detail.
In \S2 and \S3 we give the dynamical
connection underlying ``Hypothesis A.''
In \S2 we review radix expansions to an integer
 base $b \geq 2$ treated as a discrete dynamical system 
acting on the interval $[0,1]. $
The radix expansion of a real number $\theta$
is described by an orbit of a dynamical system,
the $b$-transformation $T_b (x) = b x~(\bmod ~1),$ 
 studied by Renyi~\cite{Re57} and Parry~\cite{Pa60}.
For a given number $\theta$ its $b$-expansion can be computed
from the iterates of this system
$$ x_{n+1}  = b x_{n} ~(\bmod ~ 1), $$
with initial condition $x_0 = \theta ~(\bmod ~ 1).$ The
{\em $b$-expansion} of a real number $\theta \in [0,~1]$ is
$$ \theta = \sum_{j= 1}^\infty  d_j b^{-j}, $$
in which the $j$-th digit $d_j := \lfloor bx_{j - 1} \rfloor.$
In \S3 we suppose the given real number $\theta$ is expressed as
\beql{101} 
\theta = \sum_{n = 1}^\infty  \epsilon_n b^{-n}   ,
 \eeq 
in which  $\epsilon_n$ is any sequence of real numbers with 
$\epsilon_n \to 0$ as $n \to \infty$. To this one can associate
a {\em perturbed $b$-expansion} associated to the 
{\em perturbed $b$-transformation}
\beql{102} 
y_{n+1} = b y_n + \epsilon_n ~(\bmod ~ 1), 
\eeq 
starting with an initial condition $y_0 \in [0, 1).$ 
The recurrence \eqn{102} is an infinite sequence of maps which change at
each iteration.
Associated to
this recurrence is the {\em perturbed $b$-expansion}
$$ y_0 + \theta = \sum_{j=0}^\infty \tilde{d_j} b^{-j}, $$
in which the $j$-th digit 
$\tilde{d_j} := \lfloor b y_{j} + \epsilon_{j+1} \rfloor. $
Choosing the initial condition $y_0 = 0$ gives a 
perturbed $b$-expansion of $\theta$. 
The mechanism underlying the approach of
Bailey and Crandall is that the
the $b$-expansion of $\theta$
and the perturbed $b$-expansion of $\theta$
are strongly correlated in the following sense:  
The orbit $\{ y_n~:~ n \geq 0\}$ of the perturbed $b$-transformation
with initial condition
$y_0 = 0$ asymptotically approaches the orbit $\{ x_n~:~ n \geq 0\}$
associated to the $b$-transformation with initial condition
$x_0 = \lfloor \theta \rfloor.$ (Theorem~\ref{th31a})
In particular, the orbits 
$\{ x_n~:~ n \geq 0\}$  and the 
orbit $\{ y_n~:~ n \geq 0\}$ have the same set of
limit points, and one is uniformly distributed (mod~ 1) if and only if
the other is. This implies that the
perturbed $b$-expansion of $\theta$, although
different from the $b$-expansion of $\theta$,
must have similar statistics, in various senses. This
connection is quite general, since every real number $\theta$
has representations of the form \eqn{101}.

In \S4 we consider the particular class 
of arithmetical constants treated in 
Bailey and Crandall \cite{BC00}, consisting of 
the countable set of $\theta$ given by
an expansion \eqn{101} with  $b \ge 2$ an integer and
$\epsilon_n = \frac {p(n)}{q(n)}$ where $p(x), q(x) \in \zz(x)$,
with $q(n) \neq 0$ for all $n \geq 0.$
We call such
numbers  {\em BBP-numbers,}
and call the associated
formula
$$\theta = \sum_{n = 1}^\infty \frac {p(n)}{q(n)} b^{-n},$$
a  {\em BBP-expansion to base $b$} of $\theta.$
These numbers are named after
Bailey, Borwein and Plouffe~\cite{BBP97}, who demonstrated the
usefulness of such representations 
(when $\deg(p(x)) < \deg(q(x))$)
in computing base $b$ radix
expansions of such numbers. 
We consider $BPP$-numbers having the additional restriction
$\deg(p(x)) < \deg(q(x)),$ for this condition is 
necessary and sufficient for
$\epsilon_n \to 0$ as $n \to \infty,$ 
so that the results of \S3 apply.
The number-theoretic character of $BPP$-numbers
is that they are special values (at rational points) of
functions satisfying a homogeneous linear differential
equation with integer polynomial coefficients.
We derive the result of  Bailey and Crandall that
``Hypothesis A'' implies that such $\theta$  
either are rational or are normal numbers to base $b$
(Theorem~\ref{th41}.)  This result
makes it of interest to find criteria
to determine when $BBP$-numbers are irrational, 
which we consider next.

In \S5 we  relate $BBP$-numbers 
to the theory of $G$-functions, and characterize
the subclass of $BBP$-numbers which are ``special values''
of $G$-functions.  The subject of $G$-functions
has been extensively developed in recent years,
see \cite{An89}, \cite{Bo81}, \cite{DGS94}, 
and the  special values of such functions can often be proved 
to be irrational,
see \cite{Bo81}, \cite{Ch84},\cite{Ga74}. 
We observe that 
$BBP$-numbers satisfy all but one of the properties required to be
a special value of a $G$-function
defined over the base field $\qq$.  
We then show that a $BBP$-expansion to base $b$ corresponds to a special
value of a $G$-series at $z = \frac{1}{b}$ if and only if  
the denominator polynomial 
$q(x)$ (in lowest terms) factors into linear factors over the
rationals.(Theorem~\ref{Nth43}.)
We show that if  all the roots of $q(x)$ are distinct,
then such  special values are either rational
or transcendental, using Baker's results on linear
forms in logarithms, in Theorem~\ref{Nth44}, a result
also obtained in Adikhari, Saradha, Shorey and Tijdeman
~\cite{ASST00}. We summarize
other known results on irrationality or transcendence
of special values of $G$-functions of the type in Theorem~\ref{Nth43}.
It is interesting to observe that every one of the
many examples given in \cite{BC00} is a special value of a $G$-function.
Many other interesting examples of such constants
were known earlier, for example
in 1975 D. H. Lehmer~\cite[p. 139]{Le75} observed that
$$
\sum_{n = 0}^\infty \frac{1}{(n+1)(2n + 1)(4n + 1)} = \frac {\pi}{3}.
$$

In \S6 we 
compare ``Hypothesis A'' with  a conjecture of Furstenberg
in ergodic theory, which concerning measures that are
ergodic for the
joint action of two multiplicatively independent $b$-transformations.
Both conjectures have similar conclusions, though there
seems to be no direct relation between their hypotheses.
Bailey and Crandall have found examples
of arithmetical constants which have the property of
being $BBP$-numbers to
two multiplicatively independent bases. This suggests that
one should look for further conditions under which the
two conjectures are more directly related.

In \S7 we make concluding remarks. We describe 
an empirical taxonomy of 
various classes of arithmetical constants, and 
formulate  some  alternative classes of
arithmetical constants as candidates for inclusion in
``Hypothesis A.''

%
%

\section{Radix Expansions}
\hsp
We consider radix expansions to an integer base $b \ge 2$.
They are obtained by iterating the {\em $b$-transformation}
\beql{201}
T_b (x) = bx ~ (\bmod~1) ~.
\eeq
Given a real number $x_0 \in [0,1)$, as initial condition, 
we produce the sequence of {\em remainders}
\beql{202}
x_{n+1} = bx_n ~ (\bmod~1) ~,
\eeq
with $0 \le x_{n+1} < 1$.
That is, 
\beql{203}
x_{n+1} = bx_n - d_{n+1}
\eeq
where
\beql{204}
d_{n+1} = d_{n+1} (x_0) = \lf bx_n \rf \in \{0, 1, \ldots, b-1 \}
\eeq
is called the {\em $n$-th digit of $\theta$.}
The {\em forward orbit} of $x_0$ is $\sO^+ (x_0) = \{x_n : n \ge 0 \}$ 
and we call $\{x_n\}$ the {\em remainder sequence} of the $b$-expansion.
Iterating \eqn{203} $n + 1$ times yields
\beql{204a}
x_{n+1} = b^{n + 1}  x_0 - d_{n+1} - b d_n - \cdots - b^{n} d_1 ~.
\eeq
Dividing by $b^{n+ 1}$ yields
$$x_0 = \sum_{j=1}^n d_j b^{-j} - b^{-n-1} x_{n + 1} ~.$$
Letting $n \to \infty$ yields the {\em $b$-expansion} of $x_0$,
\beql{205}
x_0 = \sum_{j=1}^\infty d_j (x_0) b^{-j} ~,
\eeq
which is valid for $0 \le x_0 < 1$.
For $\theta \in \rr$ we take $x_0 = \theta - \lf \theta \rf$ and
$d_0 ( \theta ) = \lf \theta \rf \in \zz$, thus obtaining the representation
\beql{205a}
\theta = d_0 ( \theta ) + \sum_{j=1}^\infty d_j (\theta ) b^{-j} ~,
\eeq
which is called the {\em b-expansion} of $theta.$
Note that \eqn{204a} gives
\beql{206}
x_n \equiv b^n x_0  ~(\bmod~1) \equiv  b^n \theta ~(\bmod~1) ~
\eeq
in this case.

The following property than an  initial condition $\theta$ 
may have concerns the topological dynamics of
the $b$-transformation for its iterates.

\begin{defi}\label{de21}
{\rm A real number $\theta \in [0,1)$ is {\em digit-dense to base} $b$
if, for every $m \ge 1$, every legal digit sequence of digits of
length $m$  occurs 
at least once as consecutive digits in the $b$-expansion}
$\theta = \sum\limits_{n=1}^\infty d_n (\theta ) \beta^{-n}$.
\end{defi}

The following property that an  initial
condition $\theta$ may have concerns  the metric dynamics of
the $b$-transformation for $\theta.$ It is well known
that  the  $b$-transformation 
$T_b$ has the uniform measure
(Lebesgue measure) on $[0,1]$ as its unique absolutely 
continuous invariant measure. 

\begin{defi}\label{de22}
{\rm 
A real number $\theta \in [0,1)$ is {\em normal to base} $b$ if for every
$m \ge 1$ every digit sequence
$d_{1} d_{2} \cdots d_{m} \in \{0,1, \ldots, d-1\}^m$ occurs with
limiting frequency $b^{-m}$, as given by the invariant measure $\mu_{Leb}$.}
\end{defi}

Recall that
$$\mu_{Leb} ( \{ x_0 : d_1 (x_0) \cdots d_n (x_0) = 
d_{1} d_{2} \cdots d_{m} \} ) = b^{-m} ~,
$$
where $\mu_{Leb}(S)$ denotes the Lebesgue measure of $S$.
It is well known that, for each $b \ge 2$, the set of 
$\theta \in [0,1]$ that is normal to base $b$ has full Lebesgue measure.

The properties of the digit expansion $\{ d_n ( \theta ) : n \ge 1 \}$
can be extracted from the remainder sequence $\{x_n\}$.
The following result is well known.

\begin{theorem}\label{th21}
Given an integer base $b \ge 2$ and a real number $\theta \in [0,1]$.
\begin{itemize}
\item[(1)]
$\theta$ is digit-dense to base $b$ if and only if its remainder sequence
$\{x_n (\theta ) : n \ge 1 \}$ to base $b$ is dense in $[0,1]$.
\item[(2)]
$\theta$ is normal to base $b$ if and only if its remainder sequence
$\{x_n : n \ge 1\}$ to base $b$ is uniformly distributed in $[0,1]$.
\item[(3)]
$\theta$ has an eventually periodic $b$-expansion.
if and only if its remainder sequence $\{x_n: n \ge1 \}$ to base $b$ has
finitely many limit points.
This condition holds if and only if $\{x_n : n \ge 1 \}$ 
eventually enters 
a periodic orbit of the $b$-transformation, i.e. $x_m = x_{m+p}$ 
for some $m$, $p \ge 1$. These equivalent conditions hold if and only
$\theta$ is rational.
\end{itemize}
\end{theorem}

\paragraph{Proof.}

(1). The set $I(d_{1} d_{2} \cdots d_{m} ) := \{\theta \in [0,1]:
d_1 ( \theta) \cdots d_m (\theta ) = d_{1} \cdots d_{m} \}$
is a half-open interval $[a,a + b^{-m})$ of length $b^{-m}$, and the $b^m$ 
intervals partition $[0,1]$.
Digit-denseness implies there exists some 
$x_k \in I(d_{1} \cdots d_{m} )$.
This holds for all $m \ge 1$ and generates a dense set of points.

(2). If $\theta$ has $\{x_n : n \ge 1 \}$ is uniformly distributed 
$(\bmod~1)$ then the correct
frequency of points occurs in each interval $I(d_{1} \cdots d_{m} )$, 
and this proves normality of $\theta$.
For the converse, one uses the fact that $I(d_{1} \cdots d_{m} )$ 
is a basis for the Borel sets in $[0,1)$.

(3). The key point to check is that if the limit set of 
$\{x_n : n \ge 1\}$ is finite, then this finite set forms 
a single periodic orbit of the $b$-transformation, 
and some $x_n$ lies in this orbit.
We omit details, cf.
 Bailey and Crandall~\cite[Theorem 2.8]{BC00}.~~~$\bsq$

\paragraph{Remark.} Most of the results above generalize to the
$\beta$-transformation $T_{\beta}(x) = \beta x~ ( \bmod~1)$
for a fixed real $\beta > 1;$ these maps were 
studied by Parry~\cite{Pa60}.
Associated to this map is the notion of
a {\em $\beta$-expansion} for any real number $\theta$, in
which the allowed digits are $\{ 0, 1, 2, \ldots, \lfloor \beta \rfloor \}.$
Not all digit sequences are allowed in $\beta$-expansions, but
the set of allowed digit sequences 
was characterized by Parry~\cite{Pa60}, 
see Flatto et al. \cite{FLP94} for other references.
One defines  a number $\theta$ to be {\em digit-dense
to base $\beta$} if every allowable finite digit sequence
occurs in its $\beta$-expansion..
There is a unique absolutely continuous invariant 
measure $d\mu$ of total mass one
for the $\beta$-transformation,
and one defines  a number $\theta$ to be
{\em normal to base $\beta$} if every finite block of digits
occurs in its $\beta$-expansion with the limiting frequency
prescribed by this invariant measure.
With these conventions, Theorem~\ref{th21} remains valid
for a general base $\beta$, 
except that Theorem~\ref{th21}(3)
must be taken only as characterizing 
 eventually periodic orbits of the
$\beta$-transformation. That is, the final assertion 
in (3) that $\theta$ is
rational must be dropped; it does not hold for general $\beta.$ 
For results relating
normality of
numbers in different  real bases $\beta$, see 
Brown, Moran and Pollington~\cite{BMP97}.

%
%

\section{Perturbed Radix Expansions}
\hsp
Let $b \ge 2$ be an integer, and let $\{\ep_n: n \ge 1 \}$ 
be an arbitrary sequence of real numbers satisfying
\beql{N301}
\lim_{n \to \infty} \ep_n =0 ~.
\eeq
Set
\beql{N302}
\theta = \theta (b, \{\ep_n \}) : = \sum_{n=1}^\infty \ep_n b^{-n} ~.
\eeq
We can study the real number $\theta$ using a 
{\em perturbed $b$-expansion} associated to the sequence $\{\ep_n\}$.

The {\em perturbed $b$-transformation} on $[0,1)$ is the recurrence
\beql{N303}
y_{n+1} = by_n + \ep_{n+1} \quad (\bmod~1) ~,
\eeq
with $0 \le y_{n+1} < 1$ and
with given initial condition $y_0$. 
That is, 
\beql{N304}
y_{n+1} = b y_n + \ep_{n+1} - \tilde{d}_{n+1},
\eeq
where
\beql{N305}
\tilde{d}_{n+1} = \lfloor by_n + \ep_{n+1} \rfloor \in \zz ~
\eeq
is the {\em $(n + 1)$-st digit} of the expansion.
The {\em digit sequence} $\tilde{d}_n = d_n (y_0)$ 
and {\em remainder sequence} 
$\{ y_n: n \geq 0 \}$
depend on the initial condition $y_0$. 
Since $\ep_n \to 0$, for all sufficiently large $n$, 
one has $\tilde{d}_n \in \{-1, 0,1, \ldots, b-1, b\}$.
Now \eqn{N304} iterated $n + 1$ times yields
\beql{N306}
y_{n+1} = \ep_{n+1} + b \ep_n + \cdots +
b^n \ep_1 + b^{n+1} y_0 -
\sum_{j=0}^n \tilde{d}_{n - j} b^j ~.
\eeq
Dividing by $b^{n+1}$ yields
\beql{N307}
\sum_{j=1}^n \tilde{d}_j b^{-j} =
\sum_{j=1}^n \ep_j b^{-j} +
(y_0 - b^{-n-1} y_{n+1} )~.
\eeq
Letting $n \to \infty$ yields the {\em perturbed $b$-expansion}
\beql{N308}
y_0 + \theta = \sum_{j=1}^\infty \tilde{d}_j (y_0 ) b^{-j} ~,
\eeq
valid for $0 \le y_0 < 1$.
We write $y_n = y_n (y_0)$ for the remainder sequence in \eqn{N304}

The {\em perturbed $b$-expansion}
$\{d_n^\ast (\theta ) : n \ge 1 \}$ for $\theta$ given by \eqn{N302} is
obtained by 
choosing the initial condition $y_0 =0$,
i.e. $d_n^\ast (\theta ) := \tilde{d}_n (0)$.
We  also have the {\em perturbed remainders }
$\{y_n^\ast (\theta ) : n \ge 1 \}$ given by $y_n^\ast (\theta ) = y_n (0)$.

The main
 observation of this section is that the remainders of the perturbed 
$b$-expansion of such $\theta$ are related to the 
remainders of their $b$-expansion. 

\begin{theorem}\label{th31a}
Let $b \ge 2$ be an integer and let 
$\theta := \sum_{n=1}^\infty \ep_n b^{-n}$, where $\ep_n$ are 
real numbers with $\ep_n \to 0$ as $n \to \infty$.
Let $\{y_n^\ast (\theta ) : n \ge 1 \}$ denote the associated 
perturbed remainder sequence of $\theta$, and 
$\{x_n (\theta ) : n \ge 1 \}$ 
denote the remainder sequence of its $b$-expansion. If
\beql{N309}
t_n := \sum_{j=1}^\infty \ep_{n+j} b^{-j} ~,
\eeq
then
\beql{N310}
x_n ( \theta) = y_n^\ast ( \theta ) + t_n \quad (\bmod~1) ~.
\eeq
The orbits $\{x_n (\theta ) : n \ge 1 \}$ and 
$\{y_n^\ast (\theta ) : n \ge 1 \}$ asymptotically approach 
each other on the torus $\TT = \rr/ \zz$  as $n \to \infty$.
\end{theorem}

\paragraph{Proof.}
Since $y_0 =0$, formula \eqn{N306} gives
\beql{N311}
y_{n+1} = \sum_{j=1}^n b^{n-j} \ep_j \quad (\bmod~1) ~.
\eeq
Now
$$b^n \theta = \sum_{y=1}^\infty b^{n-j} \ep_j = \sum_{j=1}^n
b^{n-j} \ep_j + t_n ~.
$$
Thus
\beql{N312}
b^n \theta = y_n + t_n \quad (\bmod~1)
~.
\eeq
For the $b$-expansion, \eqn{206} gives $b_n \theta \equiv x_n$ $(\bmod~1)$,
and combining this with \eqn{N312} yields \eqn{N310}.

Since $\ep_n \to 0$ as $n \to \infty$, we have $t_n \to 0$ as $n \to \infty$.
Thus 
$|x_n (\theta ) - y_n^\ast (\theta ) | \to 0 \quad\mbox{on}\quad \TT ~$
 as $n \to \infty.$
Note that on $\TT = \rr/\zz$ the points $\ep$ and $1-\ep$ 
are close.~~~$\bsq$

\begin{lemma}\label{le32}
Let $\{x_n : n \ge 1\}$ and $\{y_n : n \ge 1\}$ be any two sequences
in  $[0,1]$ 
with $x_n = y_n + \delta_n$ $(\bmod~1)$ with 
$\delta_n \to 0$ as $n \to \infty$.
\begin{itemize}
\item[(1)]
The sequences $\{x_n : n \ge 1 \}$ and $\{ y_n : n \ge 1 \}$ have the same
sets of limit points, provided the endpoints 0 and 1 are identified. 
\item[(2)]
The sequence $\{x_n : n \ge 1\}$ is uniformly distributed $(\bmod~ 1)$ 
if and only if $\{y_n : n \ge 1\}$ is uniformly distributed $(\bmod~1)$.
\end{itemize}
\end{lemma}

\paragraph{Proof.}

(1) This is clear since $x_{n_j} \to \psi$ implies $y_{n_j} \to \psi$ and
vice-versa, except at the endpoints $\psi =0$ or 1, which, by
convention, we identify as the same point.

(2) This is well known, see Kuipers and Niederreiter~
\cite[Theorem 1.2, p. 3]{KN74}.~~~$\bsq$ \\

One can compare the $b$-expansion $\{d_n (\theta ) : n \ge 1 \}$ and 
the perturbed
$b$-expansion $\{ d_n^\ast (\theta ) : n \ge 1 \}$ of such $\theta$.
We have
\begin{eqnarray*}
d_n (\theta ) & = &
\lfloor bx_{n-1} \rfloor \\
d_n^\ast (\theta ) & = & \lfloor b_{n-1} + \ep_n \rfloor =
\lfloor b(x_{n-1} + t_{n-1} (\bmod~1 )) + \ep_n \rfloor ~.
\end{eqnarray*}
Since $t_n \to 0$ and $\ep_n \to 0$ as $n \to \infty$, one expects 
that ``most'' digit values of the two expansions will agree\footnote
{This is an unproved heuristic
statement.  It is an open problem to prove that a 
natural density one proportion of all $n$ have 
$d_n (\theta ) = d_n^\ast (\theta)$.},
i.e. $d_n (\theta ) = d_n^\ast (\theta )$ for ``most'' 
sufficiently large values of $n$.
However there is still room for there to be infinitely many $n$ where
$d_n (\theta ) \neq d_n^\ast (\theta )$.

We next consider perturbed $b$-expansions having a finite
number of limit points, and show that they correspond to
rational $\theta$.

\begin{theorem}\label{th31}
Let $b \ge 2$ be an integer and let 
$\theta = \sum_{n=1}^\infty \ep_n b^{-n}$ with $\ep_n$ a sequence 
of real numbers with $\ep_n \to 0$ as $n \to \infty$.
The following conditions are equivalent.
\begin{itemize}
\item[(i)]
$\theta \in \qq$.
\item[(ii)]
The remainders $\{y_n^\ast ( \theta ) : n \ge 1 \}$ of the 
perturbed $b$-expansion of $\theta$ have finitely many 
limit points in $[0,1]$.
\item[(iii)]
The orbit $\{ y_n^\ast (\theta ): n \ge 1 \}$ of the 
perturbed $b$-transformation asymptotically approaches a 
periodic orbit $\{x_k : 0 \le k \le p \}$ of the 
$b$-transformation, with $T_b (x_k) = x_{k+1}$ and $T_b (x_p) = x_0$ and for
$0 \le j \le p-1$, such that
\beql{N313}
y_n = x_{j} + \delta_n ~(\bmod~1) ~ 
\quad\mbox{if}\quad n \equiv j ~(\bmod~p)
\eeq
with $\delta_n \to 0$ as $n \to \infty$.
\end{itemize}
\end{theorem}

\paragraph{Proof.}
(i)$\Rightarrow$(ii).  By Theorem \ref{th21} if 
$\theta \in \qq$ the remainders
$\{x_n (\theta ) : n \ge 1 \}$ of the $b$-transformation 
have finitely many limit points.
By Theorem \ref{th31a} and Lemma \ref{le32}
we conclude that $\{y_n^\ast (\theta ): n \ge 1\}$ has 
the same set of limit points.

\noindent
(ii)$\Rightarrow$(iii).  By Theorem \ref{th31a} and Lemma \ref{le32} 
the limit points of
$\{y_n^\ast (\theta ): n \ge 1 \}$ are the same as 
$\{x_n (\theta ) : n \ge 1 \}$.
By Theorem \ref{th21} such limit points must form 
a periodic orbit of the $b$-transformation.

\noindent
(iii)$\Rightarrow$(i).
The values $\{y_n^\ast (\theta ) : n \ge 1\}$ have limit points 
the periodic orbit
$\{x_j : 1 \le j \le n \}$ of $T_b$.
By Theorem \ref{th21},
it follows that $\theta \in \qq$.~~~$\bsq$

\paragraph{Remarks.}

(1). Any real number $\theta$
has some perturbed $b$-expansion satisfying the hypotheses
of Theorem~\ref{th31a}, so in a sense these expansions
are completely general.
It follows from  Theorem~\ref{th31} 
that Hypothesis A cannot be valid for all
such $\theta$, since there exist irrational $\theta$
that are not normal numbers.

(2). The rationality criterion of Theorem \ref{th31} is not
directly testable computationally, 
unless all $\ep_n =0$ for $n \ge  n_0$; the latter case
essentially is the same as that of a $b$-transformation.
When infinitely many $\epsilon_n \neq 0$,
then the points $\{y_n^\ast (\theta ) : n \ge 1 \}$ 
stay outside the periodic orbit for infinitely many values of $n$, 
and the role of the $\{\ep_n\}$ is to compensate for the 
expanding nature of the map
$T(x) = bx$ $(\bmod ~1)$ by providing negative feedback to 
push the iterates closer and closer to the periodic orbit.

(3). Theorem~\ref{th31a} does not extend to
$\beta$-expansions for non-integer $\beta.$ 
One can consider
\beql{N314}
\theta = \theta (\beta, \{\ep_n \}) : = \sum_{n=1}^\infty \ep_n \beta^{-n} ~.
\eeq
and define an associated perturbed $\beta$-transformation in
the obvious way. However when $\beta$ is not an integer
the analogue of Theorem~\ref{th31a}
fails to hold, since \eqn{N312} is no longer valid.
In particular, Theorem~\ref{th31a} does not extend to rational
$\beta = \frac {b}{a} > 1,$ with $a > 1.$

%
%
\section{BBP-Numbers and Hypothesis A}

We consider expansions of the following special form.

\begin{defi}\label{de41}
{\rm A {\em BBP-number to base $b$} is a real number $\theta$ 
with a representation
\beql{N401}
\theta = \sum_{n=1}^\infty \frac{p(n)}{q(n)} b^{-n} ~,
\eeq
in which $b \ge 2$ is an integer and $p(x)$, $q(x) \in \zz [x]$ 
are relatively prime polynomials,
with $q(n) \neq 0$ for each $n \in \zz_{\ge 0}$.
We call \eqn{N401} a {\em BBP-expansion to base $b$}.
}
\end{defi}

Bailey, Borwein and Plouffe \cite[p. 904]{BBP97} introduced 
this class of numbers, proving that for them
the $d$-th digit is computable\footnote
{Bailey, Borwein and Plouffe
 use the convention that
``computing the d-th digit'' means computing
is an approximation to $ b^d \theta ~(\bmod~ 1)$ 
that is guaranteed to be
within a specified distance to it $(\bmod~ 1)$. Usually
this determines the $d$-th digit, but it may not, near the
endpoints of the digit interval.}
in time at most 
$O(d \log^{O(1)} d )$ using  space at most $O(\log^{O(1)} d)$,
which is the complexity class $SC^\ast$, a
subclass of $SC$, see \cite[p. 127]{J90}.

We mainly consider $BBP$-numbers that satisfy the extra condition
\beql{N402}
\deg (q(x)) > \deg (p(x)) ~.
\eeq
This condition guarantees that 
$\ep_n = \frac{p(n)}{q(n)}  \to 0$ 
as $n \to \infty$, which allows Theorem~\ref{th31a} to be applicable.
We now  formulate two hypotheses, whose conclusions are in
terms of topological
dynamics and metric dynamics, respectively. The second of 
these is
``Hypothesis A'' of Bailey and Crandall~\cite{BC00}.

\noindent \paragraph{\bf Weak Dichotomy Hypothesis.} 
{\em Let there be given  a perturbed $b$-transformation with 
$\ep_n = \frac{p(n)}{q(n)}$ with $p(x)$, $q(x) \in \zz [x]$ and
$ \deg (q(x)) > \deg (p(x))$.
Then the orbit $\{y_n : n \ge 1 \}$ for $\theta (b, \{\ep_n \})$ 
either has finitely many limit points or else is dense in $[0,1]$.} \\

\noindent \paragraph{\bf Strong Dichotomy Hypothesis} 
{\em Let there be given  a perturbed $b$-transformation with
$\ep_n = \frac{p(n)}{q(n)}$ with $p(x)$, $q(x) \in \zz [x]$ 
and $\deg (q(x)) > \deg (p(x))$. 
Then the orbit 
$\{y_n : n \ge 1 \}$ for $\theta (b, \{\ep_n \} )$ either has
finitely many limit points or is uniformly distributed on
$[0,1].$ Equivalently, in measure theoretic terms, 
the measures $\mu_N =\frac{1}{N} \sum_{k=1}^N \delta_{y_k} $
converge in the vague topology as $N \to \infty$
to a limit measure $\mu$, which is an invariant
measure for the $b$-transformation, and which is either
a measure supported on a finite set or else 
is Lebesgue measure on $[0,1].$} \\

Bailey and Crandall~\cite{BC00} essentially established  
the following result.

\begin{theorem}\label{th41}
Let $\theta$ be a BBP-number to base $b$ whose 
associated BBP expansion satisfies
\beql{N403}
\deg (q(x)) > \deg (p(x)) ~.
\eeq
Then the following conditional results hold.

(1) The Weak Dichotomy Hypothesis implies that $\theta$ is either rational 
or digit-dense to base $b$.

(2) The Strong Dichotomy Hypothesis implies that  
$\theta$ is either rational or  a normal number to base $b$.
\end{theorem}

\paragraph{Proof.}
The condition \eqn{N403} guarantees that 
$\ep_n = \frac{p(n)}{q(n)} \to 0$ as $n \to \infty$.
Thus Theorem \ref{th31a} applies to the $BBP$-number
$$\theta = \sum_{n=1}^\infty \frac{p(n)}{q(n)} b^{-n} ~.$$

(1) By the Weak Dichotomy Hypothesis the limit set of 
$\{y_n^\ast (\theta ) : n \ge 1 \}$ is dense in $[0,1]$.
Now Lemma \ref{le32} (1) implies that $b$-expansion remainders 
$\{x_n (\theta ): n \ge 1\}$ are dense in $[0,1]$.
Theorem \ref{th21} (1) then shows that $\theta$ is digit-dense.

(2) By the Strong Dichotomy Hypothesis 
$\{y_n^\ast ( \theta ) : n \ge 1 \}$ is uniformly distributed in $[0,1]$.
Now Lemma \ref{le32} (2) implies that 
$\{x_n (\theta ) : n \ge 1 \}$ is uniformly distributed in $[0,1]$.
Now $\theta$ is normal to base $b$ by Theorem \ref{th21} (2).~~~$\bsq$

Bailey, Borwein and Plouffe \cite{BBP97} and 
Bailey and Crandall \cite{BC00} give many examples of $BBP$-numbers 
satisfying \eqn{N402} where the associated real number $\theta$ is known 
to be irrational.
For example for various $b$ one can obtain $\pi$, $\log 2$, $\zeta (3)$ etc.
They also observe that $\zeta (5)$ is a $BBP$-number, 
to base $b=2^{60}$, but it remains an open problem to
decide if $\zeta (5)$ is irrational. All the examples they
give of $BBP$-numbers are actually of a special form: they
are ``special values'' of $G$-functions defined over $\qq$,
as we discuss next.

%
%

\section{Special Values of G-Functions}
\hsp

The notion of $G$-function was introduced by Siegel~\cite{Si29}
in 1929.

\begin{defi}\label{de42}
{\rm A power series
\beql{N404}
f(z) = \sum_{n=0}^\infty a_n z^n
\eeq
defines a {\em $G$-series over the base field $\qq$} if the 
following conditions hold.

(i) {\em Rational coefficients condition}.
All $a_n \in \qq$ so we may write $a_n = \frac{p_n}{q_n}$, 
with $p_n, q_n \in \zz$ with $(p_n, q_n) =1$ and  $q_n \ge 1$.

(ii) {\em Local analyticity condition}.
The power series $f(z)$ has positive radius of convergence 
$r_\infty$, and for each prime $p$ the $p$-adic function 
$f_p (z) := \sum_{n=0}^\infty a_n z^n$ viewing 
$a_n \in \qq \subseteq \qq_p$ has positive radius of convergence $r_p$
in $\cc_p$, where $\cc_p = \hat{\bar{\qq_p}}$ is the completion of
the algebraic closure of $\qq_p.$

(iii) {\em Linear differential equation  condition}. The power series  $f(z)$
formally  satisfies a 
homogeneous linear differential equation in
$D= \frac{d}{dz}$ with coefficients in the polynomial ring
$\qq [z]$.

(iv) {\em Growth condition}.
There is a constant $C< \infty$ such that
\beql{N405}
g_n := {\rm lcm} (q_1, q_2, \ldots, q_n ) < C^n
\eeq
for all $n \ge 1$.
}
\end{defi}

\noindent
There is an extensive theory of $G$-functions, see  
Bombieri \cite{Bo81}, Andr\'{e} \cite{An89} and 
Dwork, Geroth and Sullivan \cite{DGS94}.
For the general definition of a $G$-function over an
algebraic number field 
$K$ see Andr\'{e} \cite[p. 14]{An89},
or Dwork et~al. \cite{DGS94}. 
$G$-functions have
an important role in arithmetic algebraic geometry, where
it is conjectured that 
$G$-functions are exactly the set of solutions over $\bar{\qq}[z]$
of a geometric differential equation over $\bar{\qq},$ as
defined in Andre~\cite[p. 2]{An89}. In any case it is known
 that the (minimal) homogeneous linear differential equation
satisfied by a $G$-series is of a  very restricted kind: it must
have regular singular points, and these must all have rational exponents,
by a result of Katz, cf. Bombieri~\cite[p. 46]{Bo81} and 
Bombieri and Sperber~\cite{BS82}. (The growth condition (iv) 
plays a crucial role in obtaining this result.)
It follows that  a $G$-series analytically continues to a 
multi-valued function on $\pp^{1} (\cc )$ minus a finite number 
of singular points, cf. \cite[p. xiv]{DGS94}.
We call this multi-valued function a {\em $G$-function}.

It is known that the set $\sG_K$ of $G$-series defined over a 
number field $K$ form a ring over $K$, under addition and multiplication, 
which is also closed under the Hadamard product
\beql{N406}
f \boxtimes g (z) = \sum_{n=0}^\infty a_n b_n z^n ~,
\eeq
see \cite[Theorem D, p. 14]{An89}.

\begin{defi}\label{de43}
{\rm 
A {\em special value of a $G$-function} defined over $K$ is a 
value $f(b)$, where $b \in K$, which is obtained by 
analytically continuation  along some path from 0 to $b$ that
avoids singular points.}
\end{defi}

Siegel~\cite{Si29} 
introduced $G$-functions and observed that irrationality
results could be proved for their ``special values'',
but did not give any details.
Bombieri~\cite{Bo81} developed the theory of $G$-functions
and gave explicit
irrationality criteria in specific cases (his Theorem 6)
for points close to the center of the circle of convergence
of the $G$-series,  as
a by-product of very general results.

It is easy to show that each $BBP$-number is a special value of 
a power series on $\qq$ that satisfies conditions (i)--(iii) 
of a $G$-series. They do not always satisfy
the growth condition (iv), however, and in a subsequent result we  give 
necessary and sufficient  conditions
for the condition (iv) to hold.

\begin{theorem}\label{Nth42}
Let $R(x) = \frac{p(x)}{q(x)} \in \qq (x)$ with $p(x)$, 
$q(x) \in \qq [x]$ with
$(p(x), q(x)) =1$ and with
$q(n) \neq 0$ for all $n \ge 0$, and set
\beql{N407}
f(z) = \sum_{n=0}^\infty \frac{p(n)}{q(n)} z^n ~.
\eeq
Let $f_p (z)$ be the $p$-adic power series obtained by interpreting
$\frac{p(n)}{q(n)} \in \qq \subseteq \qq_p$.
Then the  power series $f(z)$ satisfies a 
homogeneous linear differential equation in
$\frac{d}{dz}$ with coefficients in $\qq [z]$, and $f(z)$ 
has positive radius of convergence in $\cc$ and $f_p (z)$ has 
a positive radius of convergence in 
$\cc_p = \widehat{\bar{\qq}_p}$ for all primes $p$.
\end{theorem}

\noindent\paragraph{Proof.}
For the first assertion, let $p(x) = \sum_{j=0}^l a_j x^j$ and 
$q(x) = \sum_{j=0}^m b_j x^j.$
Then the operator
\beql{N407a}
D := \frac{d^{l + 1}}{dz^{l + 1}} (1 - z)^{l + 1} 
(\sum_{j = 0}^m b_j ( z \frac{d}{dz})^j)
\in \qq[z, \frac{d}{dz}]
\eeq
has the property that 
\beql{N407b}
Df(z) = 0.
\eeq
Indeed one has
$$ q( z \frac{d}{dz}) f(z) = \sum_{n = 0}^\infty p(n) z^n = 
\sum_{j = 0}^l a_j^\prime (\frac{1}{1-z})^{j+1},$$
where $a_j^\prime$ are defined by the polynomial identity
$$\sum_{j=0}^l a_j x^j = \sum_{j=0}^l a_j^\prime {\binom{x}{j}}.$$ 
Multiplying this rational function by $(1-z)^{l + 1}$
yields a polynomial of degree $l$ in $z$, which is  annihilated
by $\frac{d^{l + 1}}{dz^{l + 1}},$ and this verifies \eqn{N407b}.

For the second assertion, the power series expansion of $f(z)$
clearly has radius of convergence $1$ in $\cc$.
It is easy to establish that the 
 the $p$-adic series $f_p (z)$ has a positive radius of convergence on
some $p$-adic disk around zero  
since $|q(n) | \le cn^d$ cannot contain more than 
$cd \log n$ factors of $p$. 
~~~$\bsq$

We now give necessary and sufficient conditions for a 
power series arising from a $BBP$-number to be a
$G$-series.

\begin{theorem}\label{Nth43}
Let $R(x) = \frac{p(x)}{q(x)} \in \qq (x)$ with $p(x)$, 
$q(x) \in \qq [x]$ with
$(p(x), q(x)) =1$ and with
$q(n) \neq 0$ for all $n \ge 0$, and set
$f(z) = \sum_{n=0}^\infty \frac{p(n)}{q(n)} z^n ~.$
Then the  power-series $f(z)$ is a $G$-series $($necessarily 
defined over $\qq)$ if and only if
$q(x)$ factors into linear factors in $\qq [x]$.
\end{theorem}

\noindent\paragraph{Proof.}
 Suppose first that $q(x)$ factors into linear factors over $\qq$, say
$$q(x) = A \prod_{j=1}^r L_j (x)$$
with $L_j (x) = l_j x + m_j$ with  $l_j, m_j \in \zz$ with        
$(l_j, m_j) =1$. To show $f(z)$ is a $G$-series, 
by Theorem~\ref{Nth42} it suffices to we check the growth condition (iv).
Now
\begin{eqnarray}\label{N409A}
{\rm lcm} (q_1, q_2, \ldots, q_n ) & \le & 
{\rm lcm} (q(1), q(2), \ldots, q(n)) \nonumber \\
& \le & |A|\prod_{j=1}^r {\rm lcm} (L_j (1), \ldots, L_j (n))
\end{eqnarray}
where $L_j (n) = l_j x + m_j$.
It is well-known that
\begin{eqnarray}\label{N410}
\log ({\rm lcm} [1,2, \ldots, m ]) & = &
\sum_{\{  p, j: p^j \le m\}} \log p \nonumber \\
& = & \sum_{n=1}^m \Lambda (n) = m+O(m)
\end{eqnarray}
by the prime number theorem. This yields
\beql{N410A}
{\rm lcm} [1,2, \ldots, m] = e^{m(1+o(1))}
\eeq
as $m \to \infty$.
This gives a bound
\begin{eqnarray*}
{\rm lcm} (L_j (1), \ldots, L_j (n)) & \le & 
 {\rm lcm}(1,2, \ldots, |l_j | n + |m_j| )\\
& \le & e^{(|l_j | n + |m_j|)(1+ o(1))}~.
\end{eqnarray*}
Substituting this in \eqn{N409A} implies condition (iv).

For the opposite direction, we will
show that if $q(x)$ does not factor into linear factors
over $\qq$ then condition (iv) does not hold.
Nagell~\cite{Na29} showed that
if $q(x) \in \zz[x]$ is an irreducible
polynomial of degree $d \ge 2$, then there is a positive constant
$c(d)$ with the property that for any $\epsilon > 0$ there is
a positive constant $ C(\epsilon)$ such that
\beql{410B}
{\rm lcm} (q(1), q(2), \ldots, q(n)) > C(\epsilon)  n^{ (c(d)- \epsilon) n}
\eeq
holds for all  $n \geq 1$. One can
prove this result with $c(d) = \frac {d  - 1} {d^2}$.
Such a lower bound applies to any denominator $q(x)$ that does
not split into linear factors over $\qq$.
To complete the argument one must bound the possible cancellation 
between the numerators $p(n)$, and denominators $q(n)$.
If $(p(x), q(x)) = 1$ over
$\zz[x]$, then
\beql{410C} 
\prod_{j = 1}^n {\rm gcd} (p(j), q(j)) \leq C^n, 
\eeq
for a finite constant $C = C(p(x), q(x)).$ This follows since
$$ {\rm gcd} (p(n), q(n)) \leq C $$
holds for all $n$, for a suitable $C$.
To see this, factor  $p(x)= \prod( x - \alpha_i)$
and $q(x)= \prod (x - \beta_j)$, with $\alpha_i \neq \beta_j$
for all $i, j$. Then one has, over the number field $K$ spanned by these
roots,
\beql{410D} 
{\rm ideal-gcd}( (n - \alpha_i), (n - \beta_j) ) ~| ~(\alpha_i - \beta_j).
\eeq
Taking a norm from $K/ \qq$ of the product of all these ideals
gives the desired  constant $C$. ~~~$\bsq$

\paragraph{Remarks.}(1) It is an interesting open question to
determine what is
the largest value of $c(d)$ allowed in \eqn{410B}. 
One can prove that it cannot be larger than $d - 1$.

(2) There are many more $G$-functions defined over
$\qq$ than those given in Theorem~\ref{Nth43}. The set of
$G$-functions defined over $\qq$  is closed under multiplication, so that
$(\log(1 - z))^2)$ is a $G$-function, but its power series
coefficients around $z = 0$ are not given by a rational function.
Also, for rational $a, b, c$ the Gaussian hypergeometric function 
$$
{{}_2F}_1(a, b, c, z)= 
\sum_{n=0}^\infty \frac{(a)_n(b)_n}{(c)_n n!} z^n,$$
is a $G$-function which is not of the above kind
for ``generic'' $a, b, c,$ see Andr\'{e}~\cite{An96}.

According to the results of \S4, the conclusion of 
``Hypothesis A'' is really a statement about irrational
$BBP$-numbers.
There is a good deal known about the irrationality or transcendence
of the special values of the $G$-series covered in Theorem~\ref{Nth43},
a topic which  we now address.

\begin{theorem}\label{Nth44}
Let $R(x) = \frac{p(x)}{q(x)} \in \qq (x)$ with $p(x)$, 
$q(x) \in \qq [x]$ with
$(p(x), q(x)) =1$ and with
$q(n) \neq 0$ for all $n \ge 0$, and set
$$
f(z) = \sum_{n=0}^\infty \frac{p(n)}{q(n)} z^n ~.
$$
If $q(z)$ factors into distinct linear factors over $\qq$,
then for each rational $r$ in the open disk of convergence
of $q(z)$ around $z=0$ the special value $f(r)$ is
either rational or transcendental. Furthermore there is an effective
algorithm to decide whether $f(r)$ is rational or transcendental.
\end{theorem}

\noindent\paragraph{Proof.}
We only sketch the details, since a similar result has been
obtained by Adhikari, Saradha, Shorey
and Tijdeman~\cite{ASST00}, see also
Tijdeman~\cite[Theorem 6]{Ti00}.

By expanding $R(x)$ in  partial fractions, 
under the hypothesis that $q(x)$ splits in linear factors
over $\qq$ 
one
obtains an expansion of the form 
$$R(x) = p_0(x) + \sum_{j=1}^s \frac {c_j}{x - r_j}, $$
in which $p_0(x) \in \qq[x]$, and each $c_j, r_j \in \qq$. 
In fact $r_j \notin \zz_{\geq 0}$, so 
all denominators $q(n) \neq 0.$ Now if $r_j = \frac {p_j}{q_j}$
then one has a decomposition,
$$ \sum_{n= 0}^\infty \frac{1}{n - r_j} z^j = p_j(z)  + 
\sum_{k = 1}^{q_j} \beta_{j,k}  \log(1 - \exp(\frac{2\pi i k}{q_j}) z),$$
in which $p_j(z)$ is a polynomial with rational coefficients,
while  the
coefficients $\beta_j$ are effectively computable algebraic numbers in
the field $\qq(exp(\frac{2\pi i}{q_j})).$
It follows from this
that one  can express  the function $f(z)$ as a finite
sum of terms of the form  $\frac {a_j}{(1 - z)^j}$
with rational coefficients  
plus a finite sum of terms of the form $- \beta_{j,k} \log( 1 - \alpha_j z)$,
with $\beta_j, \alpha_j$ effectively computable algebraic numbers.
The non-logarithmic terms all combine to give a rational
function $R_0(z)$ with coefficients in $\qq$.
Given a rational 
$r$ with $ 0 < |r| < 1$, it
 follows that $f(r)$ is a finite sum of linear forms in
logarithms with algebraic coefficients, evaluated at
algebraic points. Using  Baker's transcendence result
on linear forms in logarithms (Baker~\cite[Theorem 2.1]{Ba75}),
$f(r)$ is transcendental if and only if  the sum of all the
logarithmic terms above is nonzero.
There is also an effective decision procedure to tell whether this sum
is zero or not.
If the logarithmic terms do sum to zero, then the
 remaining rational function terms sum up
to the rational number $f(r)= R_0(r).$~~~$\bsq$

The case where $q(x)$ factors into linear factors over $\qq$
but has repeated factors is not
covered in the result above. 
This case includes
the polylogarithm ${\rm Li}_k(z) = \sum_{n = 1}^\infty \frac{z^n}{n^k}$
of order $k$,
for each $k \geq 2$. 
Various results are
known concerning the irrationality of such numbers.
For example, ${\rm Li}_k(\frac{1}{b})$
is irrational for all sufficiently large integers $b$, see
Bombieri~\cite{Bo81}. In fact it is known that the set of
numbers $1, {\rm Li}_1(\frac{p}{q}),..., {\rm Li}_n(\frac{p}{q})$,
with ${\rm Li}_1(z) = \log(1 - z)$, are
linearly independent over the rationals whenever
$|p| \geq 1$ and $|q| \geq (4n)^{n(n-1)} |p|^{n},$ 
according to Nikishin~\cite{Ni79}.
For polylogarithms
one has ${\rm Li}_{k}(1) = \zeta(k)$, also on the
boundary of the disk of convergence. It is not
known whether $\zeta(k)$ is irrational for odd $k \geq 5$,
although a very recent result of T. Rivoal~\cite{Ri00} shows
that an infinite number of $\zeta(k)$ for odd $k$ must
be irrational.

%
%

\section{Invariant Measures and Furstenberg's Conjecture}
\hsp
It is well known that for single expanding dynamical system, such as the
$b$-transformation $T_b$ , there always exist chaotic orbits 
exhibiting a wide range of pathology.
For example, there exist uncountably many
$\theta \in [0,1]$ whose $2$-transformation 
iterates $\{x_n\}$ satisfy
$$\frac{1}{25} < x_n < \frac{24}{25} \quad\mbox{for all}\quad
n \ge 0 ~,
$$
see Pollington \cite{Po79}. One can obtain ergodic invariant
measures of $T_b$ supported on the closure of suitable orbits,
which for example may form Cantor sets of measure zero.

If one considers instead
two  $b$-transformations, say  $T_{b_1}$
and $T_{b_2},$ with multiplicatively independent values,  
i.e. which  generate a non-lacunary commutative semigroup
$\sS = \langle T_{b_1}, T_{b_2} \rangle$, then the set of 
ergodic invariant measures for the whole semigroup 
is apparently of an extremely restricted form.
Furstenberg has proposed the following conjecture, 
suggested as an outgrowth of his work on 
topological dynamics, cf. Furstenberg \cite[Sect. IV]{Fu67}.
It is explicitly stated in Margulis~\cite[Conjecture 4]{Ma00}.
\paragraph{Furstenberg's Conjecture.}
{\em Let $a,b \ge 2$ be multiplicatively independent integers.
The only Borel measures on $[0,1]$ that are 
simultaneously invariant ergodic measures for $T_a (x) = ax$ 
$(\bmod ~1)$ and $T_b (x) = bx$ 
$(\bmod~1)$ are Lebesgue measure and measures supported on finite sets 
which are periodic orbits of both $T_a$ and $T_b$.
} \\

Various results concerning this conjecture appear in Rudolph \cite{Ru90}, 
Parry \cite{Pa96},
Host \cite{Ho95} and Johnson \cite{Jo92}. In particular, if there
is any exceptional invariant measure violating the
conjecture, it must have entropy zero
with respect to Lebesgue measure.

Furstenberg's conjecture involves some ingredients similar to
``Conjecture A'', and its conclusion involves a dichotomy 
similar to that in ``Conjecture A.''
This makes it natural to ask if there is any relation between
the two conjectures. At present none is known, in either direction.

One may look for $BBP$-numbers
$\theta \not\in \qq$ which have properties similar to that expressed
in the hypothesis of Furstenberg's conjecture
, i.e. which possess $BBP$-expansions
to two multiplicatively independent bases.
 It is known that there exist irrational 
$BBP$-numbers 
$\theta = \sum_{n = 1}^\infty R(n) b^{-n}$ 
which do possess $BBP$-expressions to two multiplicatively independent bases.
For example, Bailey and Crandall~\cite{BC00} observe that  
$\theta = \log 2$ has this property, on taking  
$$ b= 2 ~\qquad\mbox{and}\qquad~  R(x) = \frac {1}{x},$$
and 
$$b = 3^2 \qquad\mbox{and}\qquad R(x) = \frac{6}{2x - 1},$$
see \cite[eqn (4), and (10)]{BC00}. They
also observe that
$\theta = \pi^2$
has this property, as it possesses $BBP$-expansions to bases $b = 2$ and 
$b = 3^4$, the latter one found by 
Broadhurst~\cite[eqn. (212), p. 35]{Br98}. \\

{\em Question.} Do all $BBP$-numbers which are special values of $G$-functions 
have $BBP$-expansions in two multiplicatively independent bases? \\

To make tighter a possible connection between the two conjectures, one
can ask for which numbers does the following weaker version 
of ``Hypothesis A'' hold. \\

{\bf Invariant Measure Hypothesis}
{\em Every BBP-number to base $b$  has  $b$-transformation
iterates $\{ x_n\}$ that are asymptotically distributed according
to a limiting measure on $[0, 1].$} \\

It would be interesting to find extra hypotheses on
a class of arithmetical constants under 
which a precise connection can be established between
``Hypothesis  A'' and Furstenberg's conjecture.

%
%

\section{Concluding Remarks}

Many of the examples of  arithmetical constants 
arise as special values of $G$-functions defined
over the rationals, or at least ``special values'' of
functions satisfying linear differential equations with
polynomial coefficients in $\qq[x]$.
Based on the known results, one may empirically group these constants into
three classes, of apparantly increasing order of difficulty
of establishing irrationality or transcendence results.

(1).
special values of $G$-functions
$f(\frac{p}{q})$ defined over the rationals,
with $\frac{p}{q}$ inside the disk of convergence of the
$G$-series. 

(2). 
``singular values'' $f(1)$ of such a
$G$-function, which are values
taken  at a singular point of the associated 
(minimal order) linear differential
equation, on the boundary of the disk
of convergence of a $G$-series,
at which  the $G$-expansion converges absolutely.

(3).
``renormalized singular values,'' which are the
constant terms in an asymptotic expansion of a $G$-function
around a singular point.

In this hierarchy, an  arithmetical
 constant may occur as  more than one type. For
example, $\frac{\pi^2}{6} = \zeta(2)= {\rm Li}_2(1)$ occurs
as a number of type (2), but it is also  realized 
as a number of type (1), which falls in the
class of constants considered in this paper.
 It is a nontrivial problem to
determine what is the lowest level in the hierarchy
a given constant belongs.

Various constants of types (1) and (2) appear in
the renormalization of massive Feynman diagrams, see
Broadhurst~\cite{Br98} and Groote, K\"{o}rner and Pivovarov~\cite{Gr00},
who cite ${\rm Li}_4(\frac{1}{2})$ as such a constant.
Multiple zeta values and polylogarithms give many examples of type (2), see
Borwein et al. ~\cite{BBB97}, \cite{BBBL00}.
Many of the most interesting
arithmetical constants naturally arise
as constants of type (2) and (3). 
For examples, the values $\zeta(k) = {\rm Li}_k(1)$ appear as
constants of type (2), while
 Euler's constant  appears as a type (3)
``renormalized'' value at $z = 1$ of ${\rm Li}_1(z)$.
The problem of showing  the linear independence of
all odd zeta values $\zeta(2n + 1)$ over the rationals
has recently been of  great interest 
from connections with various conjectures
in arithmetical algebraic geometry, see
Goncharov~\cite{Go00}. 
Many other examples of type (2) and (3) constants appear in
Lehmer~\cite{Le75} and Flajolet and Salvy~\cite{FS98}.
I am not aware of any irrationality
or transcendence results proved for a constant of type (3).

One can extend the hierarchy above outside the class
of $G$-functions. E. Bombieri observes that the power series
$$ h(z) = \sum_{n = 1}^\infty \frac{1}{ n(n^2 + 1)} z^n $$
of $BBP$-type, which is not a $G$-series,
has special value at $z = 1 $ given by 
$$h(1) = \frac{1}{2} \Re{( \frac {\Gamma' (i)} {\Gamma (i)})}.$$
The value $z=1$ lies on the boundary of the 
disk of convergence of the power series for this
function, and corresponds to type (2) above. Another example is
$$ 
\sum_{n = 1}^\infty \frac {(-1)^n}{n^2 + 1} = 
\frac{2\pi}{e^\pi - e^{-\pi}}-1,
$$
see Flajolet and Salvy\cite[p. 18]{FS98}, who
give many other interesting examples.

The relevant special values of a rational power series 
for the approach of
Bailey and Crandall to apply are $z = \frac{1}{b}$
for integer $b \geq 2$, where the disk of convergence of
the associated power series has radius $1$. One observes
that the  theory of $G$-functions
provides irrationality results for rational values $z = \frac{a}{b}$,
without regard for whether $a = 1$ or not. This suggests
the following question. \\

\noindent{\em Question.} Given a rational value
$z = \frac{a}{b}$, with $1 < |a| < |b|$, is 
 there an associated dynamical system
(possibly higher dimensional) for which an analogue of
Theorem~\ref{th31} holds, relating the dynamics of
one orbit to 
the $\beta$-expansion of  
$\theta$, with $\beta = \frac{a}{b}.$? \\

At  present there seems to be no
 evidence that strongly favors a particular
class of arithmetical  constants for which ``Hypothesis A'' might
be expected to hold. The discussions of \S5 and
\S6 suggest that one might consider the following
classes.

(1). The largest class is the  set of ``special values'' of 
power series $f(z)$ defined over $\qq$ at $z = \frac{1}{b}$,
arising 
from solutions of $Df(z) = 0$ for some $D \in \sW := \qq [z, \frac{d}{dz}]$,
whose power-series coefficients $a_n \to 0$ as $n \to \infty.$
This class includes all BBP-numbers.

(2) One could restrict to the subclass of
special values $z = \frac{1}{b}$ 
of $G$-functions defined over the rationals.
However we know of no compelling reason to restrict to 
special values of $G$-functions. 

(3) The smallest class consists of a class of
arithmetical constants which satisfy extra conditions analogous
to the hypotheses of 
Furstenberg's conjecture. These consist of those constants
which are $BBP$-numbers to at least two multiplicatively independent
bases. One might add the further
restriction that they also be special values of G-functions.
As noted in \S6, this class includes $\pi^2$ and $\log 2.$

\paragraph{Acknowlegments.} The author thanks 
E. Bombieri for helpful
information  concerning $G$-functions, and for suggesting the
argument establishing \eqn{410C} in Theorem~\ref{Nth43}. 
He thanks D. H. Bailey for references, and  J. A. Reeds and
the referee for
helpful comments. Work on this paper was done in part
during a visit to the Mathematical Sciences Research Institute, Berkeley,
Sept. 2000.

\noindent{\rm AT\&T Labs--Research, Florham Park, NJ 07932-0971, USA} \\
{\em
email:} {\tt jcl@research.att.com }

\end{document}